\documentstyle[11 pt]{article} \textwidth 6.1 in \textheight 9.5 in

\title{\bf  Positive solutions of Robin problem for semilinear elliptic equations and a threshold result \\ \footnotetext{* Corresponding author. E-mail address:daiqiuyi@yahoo.com.cn(Q.Dai)}\footnotetext{This work is supported by NNSFC(Grant: No.10971061)}}
\author{ Junhui Xie,\ \ Qiuyi Dai$^*$,\ \ Huaxiang Hu \\
\small {\it Department of Mathematics, Hunan Normal University}\\
\small {\it Changsha Hunan 410081, PR China}\\
}\baselineskip 0.2in
\date{}
\begin{document}
\maketitle
\begin{center}
\begin{minipage}{130mm}
{\small {\bf Abstract}\ \ Let $\Omega$ be a bounded domain in $R^n$.
The main topics of this paper are the following Robin problem for
semilinear elliptic equation
\begin{equation}
\left\{\begin{array}{ll}\label{eq1}
 - \Delta u=u^p+f(x), &x\in\Omega,\\
u>0, &x\in\Omega,\\
\frac{\partial u}{\partial \nu}+\beta u=0,&x\in \partial \Omega,
\end{array}
\right.
\end{equation}
and its parabolic version
\begin{equation}
\left\{\begin{array}{ll}\label{eq2}
u_t-\Delta u=u^p+f(x), &(x,t)\in\Omega\times(0,T),\\
\frac{\partial u}{\partial \nu}+\beta u=0,&(x,t)\in \partial \Omega\times[0,T),\\
u(x,0)=u_0(x)\geq 0,&x\in\Omega.\\
\end{array}\right.
\end{equation}
We prove that for any given function $f(x)$ satisfies
$(\mathcal{F})$ displayed in the introduction, there exists a
positive number $\beta^\star_f$ such that problem (\ref{eq1}) has no
solution when $\beta\in(0,\beta^\star_f)$, and has at least two
solutions when $\beta>\beta^\star_f$. Moreover, among all solutions
of problem (\ref{eq1}) there is a minimal one. Concerning problem
(\ref{eq2}), we prove that the minimal solution of problem
(\ref{eq1}) is stable, whereas, any other solution of problem
(\ref{eq1}) is an initial datum threshold for the existence and
nonexistence of global solutions to problem (\ref{eq2}). For more
details, see Theorem 1.1-1.4 in the introduction.

}
\end{minipage}
\end{center}
\baselineskip 0.2in

\vskip 0.2in

\section{Introduction}

\setcounter{section}{1}

\setcounter{equation}{0}

\vskip 0.1in

Let $n\geq 2$ and $\Omega$ be a bounded domain in $R^n$. Denoting by
$\partial\Omega$ the boundary of $\Omega$, we consider the following
problem
\begin{equation}
\left\{\begin{array}{ll}\label{eq3}
 - \Delta u=u^p+f(x),&x\in\Omega,\\
u>0,&x\in\Omega,\\
\frac{\partial u}{\partial \nu}+\beta u=0,&x\in \partial \Omega,
\end{array}
\right.
\end{equation}
where $\nu$ is the outward unit vector normal to $\partial\Omega$,
$\beta\geq0$ is a parameter and $f(x)\geq0, f(x)\not\equiv0$ is a
given function in $C^1(\overline\Omega)$. \vskip 0.1 cm Problem
(\ref{eq3}) has different names depending on  the different values
of the parameter $\beta$. It is called Dirichlet if $\beta=+\infty$,
is called Neumann in the case $\beta=0$ and is called Robin provided
that $0<\beta<+\infty$. \vskip 0.1 cm It is worth pointing out that
problem (\ref{eq3}) occurs in various branches of mathematical
physics and biological models. Theoretically, the main topics in the
study of problem (\ref{eq3}) are to investigate the structure of
solution sets and the stability or instability of its solutions.
\vskip 0.1 cm In the case $\beta=0$, it is trivial to see that
problem (\ref{eq3}) has no solution. In the case $\beta=+\infty$,
there are many literature on problem (\ref{eq3}). If $f(x)\equiv0$,
then problem (\ref{eq3}) is reduced to the following problem
\begin{equation}
\left\{\begin{array}{ll}\label{eq4}
 - \Delta u=u^p,\hspace{0.5cm}
&x\in\Omega,\\
u>0,&x\in\Omega,\\
u=0,&x\in \partial \Omega.
\end{array}\right.
\end{equation}
It is well-known now that the existence of solutions to problem
(\ref{eq4}) strongly depends on the range of $p$. If $0<p<1$, then
problem (\ref{eq4}) has a unique solution, if $1<p<\frac{n+2}{n-2}$,
then problem (\ref{eq4}) has at least one solution, if
$p\geq\frac{n+2}{n-2}$ and $\Omega$ is a star shaped domain then
problem (\ref{eq4}) has no solution. Uniqueness and multiple results
have also been obtained by many authors for the case
$1<p<\frac{n+2}{n-2}$ and some special domains (see for example
\cite{ref LiYanyan,ref Cvc,ref LMF}). A celebrate result given by
Gidas and Spruck \cite{ref BJ} says that the solution set of problem
(\ref{eq4}) is compact in $L^\infty(\Omega)$ when
$1<p<\frac{n+2}{n-2}$, though it may have many solutions. In 1992,
Tarantello investigate the effect of inhomogeneous term $f(x)$ on
the structure of solution set of problem (\ref{eq4}), she considered
in \cite{ref TRA} the following problem
\begin{equation}
\left\{\begin{array}{ll}\label{eq5}
 - \Delta u=|u|^{p-1}u+f(x),\hspace{0.5cm}
&x\in\Omega,\\
u=0,&x\in \partial \Omega,
\end{array}\right.
\end{equation}
and proved that problem (\ref{eq5}) admits at least two solutions
$u_0(x)$, $u_1(x)$ $\in H_0^{1}(\Omega)$ provided that
$1<p\leq\frac{n+2}{n-2}$, $f\not\equiv 0$ and
$$\int_\Omega fudx \leq
C_{n}\|\nabla u\|_{L^{2}(\Omega)}^{\frac{n+2}{2}} \ \ \ \ \ \forall
\ u\in H_0^{1}(\Omega),$$ where
$C_{n}=\frac{4}{n-2}(\frac{n-2}{n+2})^{(n+2)/4}$. Furthermore
$u_0(x)\geq0$, $u_1(x)\geq0$ when $f(x)\geq0$. Roughly speaking, the
appearance of the inhomogeneous term $f(x)$ increase the number of
the solution. For more results on problem (\ref{eq5}) see \cite{ref
DaiG, ref DaiY, ref DengY}

\vskip 0.1in

Compared with the Dirichlet Problem, there are few results on Robin
problem. However, it is worth mentioning that Gu and Liu \cite{ref
GuLiu} studied the following problem
\begin{equation}
\left\{\begin{array}{ll}\label{eq6}
 - \Delta u=u^p,\hspace{0.5cm}
&x\in\Omega,\\
u>0,&x\in\Omega,\\
\frac{\partial u}{\partial \nu}+\beta u=0,&x\in \partial \Omega,
\end{array}\right.
\end{equation}
and proved that the solution set of problem (\ref{eq6}) is compact
in $L^\infty(\Omega)$. From the compactness of the solution set and
the fixed point theorem, they also proved that problem (\ref{eq6})
has at least one solution for any $0<\beta<+\infty$ provided that
$1<p<\frac{n+2}{n-2}$. It is also worth pointing out that some
existence, uniqueness and multiple results for problem (\ref{eq6})
with $p\in(1,\frac{n+2}{n-2}]$ can be founded in \cite{ref DaiF,ref
FuDai, ref XWang}.

\vskip 0.1 cm

In this paper, we focus our attention on the Robin problem with
sub-critical Sobolev exponent $p$. Hence, we always assume that
$0<\beta<+\infty$ and $1<p<\frac{n+2}{n-2}$ in the following
paragraph. Our purpose are two-folds, one is to investigate the
presence of the inhomogeneous term $f(x)$ how to change the
structure of solution set of problem (\ref{eq6}), the other is to
investigate the stability and instability of solutions of problem
(\ref{eq3}). Concerning the first issue, we find out that, unlike in
the Dirichlet problem, the presence of the inhomogeneous term $f(x)$
in the Robin problem can increase, as well as, decrease the number
of solutions subject to the range of the parameter $\beta$.
Concerning the second issue, we find out that the minimal solution
of (\ref{eq3}) is stable and any other solution of problem
(\ref{eq3}) is a initial data threshold for the existence and
nonexistence of global solutions to the parabolic version of problem
(\ref{eq3}).

\vskip 0.1 cm

To state our results rigorously, we introduce the following
notations and conditions imposed on $f(x)$.

\vskip 0.1 cm

Let $h(x)$ be the unique solution of the following problem
\begin{equation}
\left\{\begin{array}{ll}\label{eq7}
 - \Delta h=1,&x\in\Omega,\\
h=0,&x\in \partial \Omega.
\end{array}
\right.
\end{equation}
Setting
\begin{equation}
\left\{\begin{array}{ll}\label{eq8}
M_h=\max\limits_{x\in \Omega} h^p(x),\\
\Lambda=(\frac{1}{pM_h})^{\frac{1}{p-1}}.
\end{array}
\right.
\end{equation}
Then, a simple computation shows that
\begin{equation}\label{eq9}
\Lambda-\Lambda^{p}M_h=\frac{p-1}{p}(\frac{1}{pM_h})^{\frac{1}{p-1}}>0.
\end{equation}
For some results in this paper, the condition we impose on $f(x)$
can be read as

\vskip 0.1in

$(\mathcal{F})$ \ \ \ \ \ \ \ \  \ \ \ \ \   $f(x)\not\equiv 0 \ \ \ \mbox{and} \ \ \ 0 \leq
f(x)< \frac{p-1}{p}(\frac{1}{pM_h})^{\frac{1}{p-1}}$.

\vskip 0.1in

The first result of this paper can be stated as

\vskip 0.1in

{\bf Theorem 1.1}. For any given $f(x)\geq0$ and $f(x)\not\equiv0$.
There exists a positive number $\beta_f$ such that problem (\ref{eq3}) has
no solution for any $\beta\in(0,\beta_f)$. \vskip 0.1 cm

\vskip 0.1in

The second result is about the existence of minimal solution and
multiple result of problem (\ref{eq3}). It can be presented as

\vskip 0.1in

{\bf Theorem 1.2}. Assume that $f(x)$ satisfies $(\mathcal{F})$.
Then there exists a positive number $\beta_f^\star$ such that

\vskip 0.1 cm

(i)  Problem (\ref{eq3}) has no solution when
$\beta\in(0,\beta_f^{\star})$.

\vskip 0.1 cm

(ii) If $\beta\geq\beta_f^\star$, then problem (\ref{eq3}) has a minimal
solution $U_\beta(x)$ in the sense that for any solution
$u_\beta(x)$ of problem (\ref{eq3}) we have $U_\beta(x)\leq u_\beta(x)$.
Moreover $U_\beta$ is decreasing with respect to the parameter
$\beta$.

\vskip 0.1 cm

(iii) If $\beta>\beta_f^\star$, then problem (\ref{eq3}) has at least two
solutions.

\vskip 0.1in

{\bf Remark 1}. Let $\Omega_R=\{x\in R^n \ \|\ R<|x|<1\}$, and
$U_\beta(x)$ be the minimal solution of problem (\ref{eq3}) with
$\Omega=\Omega_R$. Setting $u(x)=v(x)+U_\beta(x)$, then $v(x)$
satisfies
\begin{equation}
\left\{\begin{array}{ll}\label{eq10}
 - \Delta v=(U_\beta(x)+v)^p-U_\beta^p(x),& x\in\Omega,\\
\frac{\partial v}{\partial \nu}+\beta v=0,&x\in \partial \Omega.
\end{array}
\right.
\end{equation}
As in \cite{ref FuDai}, we can prove that for any fixed integer $k$, there
exist at least $k$ positive solutions to problem (\ref{eq10}) provided that
$R$ is close enough to $1$ and $\beta$ is sufficiently large. As a
result, problem (\ref{eq3}) has at least $k+1$ solutions in this case.

\vskip 0.1in

Next, we consider the parabolic version of problem (\ref{eq3}), that is,
we consider
\begin{equation}
\left\{\begin{array}{ll}\label{eq11}
u_t-\Delta u=u^p+f(x),&(x,t)\in\Omega\times(0,T),\\
\frac{\partial u}{\partial \nu}+\beta u=0,&(x,t)\in \partial \Omega\times[0,T),\\
u(x,0)=u_0(x)\geq 0,&x\in\Omega,\\
\end{array}\right.
\end{equation}
when $\beta=+\infty$ and $f(x)\equiv0$, problem (\ref{eq11}) is reduced to
the following well studied problem
\begin{equation}\label{eq12}
\left\{\begin{array}{ll}\label{eq12}
u_t-\Delta u=u^p,&(x,t)\in\Omega\times(0,T),\\
u=0,&(x,t)\in \partial\Omega\times[0,T),\\
u(x,0)=u_0(x)\geq 0,&x\in\Omega.\\
\end{array}
\right.
\end{equation}
The a priori bound and decaying estimate of global solutions of
(\ref{eq12}) can be found in \cite{ref Giga} and \cite{ref PQS1, ref
PQS2}. A threshold result of (\ref{eq12}) under the condition that
(\ref{eq4}) has an unique solution was displayed in \cite{ref
Cazenave}, similar result can also be found in \cite{ref MLSN}.
Unfortunately, the uniqueness result is, in general, not true for
problem (\ref{eq4})(see \cite{ref LiYanyan,ref Cvc,ref Lin1,ref
Lin2}). It is worth pointing out that in the case $\Omega=R^n$ many
interesting results (including threshold result) of problem
(\ref{eq12}) can be found in \cite{ref DengLY, ref XWang,ref
GuiNiWang1,ref GuiNiWang2}. In this paper, we try to generalize the
threshold result of problem (\ref{eq12}) to problem (\ref{eq11}).
The main feather of our results is that we need not put uniqueness
restriction on the steady
state.\\

{\bf Theorem 1.3}. For any given $f(x)\geq0$ and $f(x)\not\equiv0$,
let $\beta_f$ be the positive number obtained in Theorem 1.1. Then
for any initial data $u_0(x)\geq0$, the local solution $u(x,t,u_0)$
of problem (\ref{eq11}) always blows up in finite time, that is, there
exists $0<T<+\infty$ such that
$$\displaystyle\lim_{t\rightarrow T^-}\displaystyle\max_{x\in\Omega}u(x,t,u_0)=+\infty.$$

{\bf Theorem 1.4}. Assume that $f(x)$ satisfies $(\mathcal{F})$ and
$1<p<\frac{n+2}{n-2}$. Let $\beta_f^\star$ denote the positive number
obtained in Theorem 1.2  and  $U_\beta$ be the minimal solution of
problem (\ref{eq3}). If $u_\beta(x)$ is an arbitrary solution of problem
(\ref{eq3}) which is distinct to $U_\beta$, then we have

\vskip 0.1 cm

(i) For any initial data $u_0(x)$, problem (\ref{eq11}) has no global
solution when $0<\beta<\beta_f^\star$.

\vskip 0.1 cm

(ii) If $\beta>\beta_f^\star$, $0\leq u_0(x)\leq u_\beta(x)$ and
$u_0(x)\not\equiv u_\beta(x)$, then problem (\ref{eq11}) has a
global solution $u(x,t,u_0)$. Moreover,
$$\displaystyle\lim_{t\rightarrow+\infty}u(x,t,u_0)=U_\beta.$$

\vskip 0.1 cm

(iii) If $\beta>\beta_f^\star$, $u_0(x)\geq u_\beta(x)$ and
$u_0(x)\not\equiv u_\beta(x)$, then the solution of problem
(\ref{eq11}) always blow up in finite time.

\vskip 0.1in

{\bf Remark 2}. Roughly speaking. Theorem 1.4 implies that the
minimal solution of problem (\ref{eq3}) is stable, while any other
solution of (\ref{eq3}) is an initial datum threshold for the
existence and nonexistence of global solution to problem
(\ref{eq11}).

\vskip 0.1in

{\bf Remark 3}. In the case $f(x)\equiv0$, problem (\ref{eq11}) becomes
\begin{equation}
\left\{\begin{array}{ll}\label{eq13}
 u_t-\Delta u=u^p,&(x,t)\in\Omega\times(0,T),\\
\frac{\partial u}{\partial \nu}+\beta u=0,&(x,t)\in \partial \Omega\times[0,T),\\
u(x,0)=u_0(x)\geq 0,&x\in\Omega.\\
\end{array}
\right.
\end{equation}
Concerning problem (\ref{eq13}), we have\\

{\bf Theorem 1.5}.  Assume that $1<p<\frac{n+2}{n-2}$. Let $U(x)$ be
an arbitrary solution of problem (\ref{eq6}). Then we have

\vskip 0.1 cm

(i) If $0\leq u_0(x)\leq U(x)$ and $u_0(x)\not\equiv U(x)$, then
problem (\ref{eq13}) has a global solution $u(x,t;u_0)$, moreover,
$$\displaystyle\lim_{t\rightarrow+\infty}u(x,t;u_0)=0.$$

\vskip 0.1 cm

(ii) If $u_0(x)\geq U(x)$ and $u_0(x)\not\equiv U(x)$, then the
solution $u(x,t;u_0)$ of problem (\ref{eq13}) always blow up in finite
time.

Unlike in \cite{ref Cazenave}, our threshold result Theorem 1.5 need
not put the uniqueness restriction on problem (\ref{eq6}).

\section{The proof of Theorem 1.1}

\setcounter{section}{2}

\setcounter{equation}{0}

\vskip 0.1in
This section is devoted to prove Theorem 1.1. To this end, we first prove some lemmas
needed in the proof of Theorem 1.1 as follows.\\

{\bf Lemma 2.1}.(\cite{ref BJ}) If $u(x)$ is a solution of problem
\begin{equation}
\left\{\begin{array}{ll}\label{eq14}
 - \Delta u=u^p,
&x\in R^n,\\
u\geq 0,&x\in R^n,\\
\end{array}
\right.
\end{equation}
with $1<p<\frac{n+2}{n-2}$, then $u(x)\equiv0$.

\vskip 0.1in

{\bf Lemma 2.2}. For $\beta\in(0,1)$ small enough, there exists a
positive constant $C$ independent of $\beta$ such that any solution
$u=u_\beta$ of problem (\ref{eq3})  satisfies
$\|u_\beta\|_{C^{2,\alpha}(\overline\Omega)}\leq C$.

\vskip 0.1in

{\bf Proof}. By elliptic theory, it is enough to show that
$\|u_\beta\|_{{L^\infty}(\overline\Omega)}\leq C$. We argue by
contradiction. Suppose that the conclusion is not true. Then there
would exist a sequence $0<\beta_j<\beta_f$ with
$\beta_j\rightarrow0$ as $j\rightarrow\infty$, a corresponding
sequence of solutions $u_j=u_{\beta_j}$ of problem (\ref{eq3}), and a
sequence of points $x_j$ in $\Omega$ such that
\begin{equation}
M_j=\|u_j\|_{{L^\infty}(\overline\Omega)}=u_j(x_j)\rightarrow\infty\ \ \ \ \ \ \mbox{as}\ \ \ \ j\rightarrow\infty.
\end{equation}
\vskip 0.1 cm

Let us consider the auxiliary function
$$v_j(y)=M_j^{-1}u_j(x_j+M_j^{-{\frac{p-1}{2}}}y),$$  which is defined
on $\Omega_j=M_j^{\frac{p-1}{2}}(\Omega-x_j)$.

It is easy to verify that $v_j(y)$ satisfies
\begin{equation}
\left\{\begin{array}{ll}\label{eq15}
 - \Delta v_j=v_j^p+M_j^{-p}f(x_j+M_j^{-{\frac{p-1}{2}}}y),\hspace{0.5cm}
&y\in\Omega_j,\\
v_j(0)=1,\\
\frac{\partial v_j(y)}{\partial
\nu}+\beta_jM_j^{-{\frac{p-1}{2}}}v_j(y)=0,&y\in
\partial \Omega_j,\\
0\leq v_j(y)\leq 1,&y\in\Omega_j. \\
\end{array}
\right.
\end{equation}

\vskip 0.1 cm

We denote by D either the whole space $R^n$ or the half space
$R_+^{n}=\{x=(x_1,x_2,\cdots,x_{n-1},x_n)\in R^n\ \|\ x_n>0\ \}$.
For any compact domain $K\subset D$, there exists a positive integer
$J$ such that $K\subset\Omega_j$ when $j>J$ due to
$\Omega_j\rightarrow D$ as $j\rightarrow +\infty$. Noting that
$0\leq v_j(y)\leq1$ and $f(x)\in C^1(\overline\Omega)$, it follows
from the standard elliptic estimate that  there exists a positive
constant $C$ independent of $j$ such that
$\|v_j\|_{C^{2,\alpha}(K)}\leq C$ with $\alpha\in(0,1)$ for all
$j>J$. By diagonal method, up to a subsequence, we may assume that
$v_j$ converges uniformly to a $C^2$ function $v$ on any compact
subset of D.

\vskip 0.1 cm

Since $M_j^{-p}f(x_j+M_j^{-{\frac{p-1}{2}}}y)\rightarrow 0$ as
$j\rightarrow +\infty$, by sending $j$ to $+\infty$ in problem (\ref{eq15}),
we see that $v$ satisfies
\begin{equation}
\left\{\begin{array}{ll}\label{eq16}
 -\Delta v=v^p,\hspace{0.5cm}&y\in D,\\
v\geq 0,\\
\frac{\partial v}{\partial \nu}=0,&y\in
\partial D\ \  if D=R_+^{n},\\
v(0)=1.\\
\end{array}
\right.
\end{equation}

\vskip 0.1 cm

If $D=R^{n}$, then $v(x)$ is a solution of
\begin{equation}
\left\{\begin{array}{ll}\label{eq17}
 -\Delta v=v^p,\hspace{0.5cm}&y\in R^{n},\\
v\geq 0,&y\in R^{n},\\
v(0)=1,\\
\end{array}
\right.
\end{equation}
this contradicts to Lemma 2.1.

\vskip 0.1 cm

If $D=R_+^{n}$, we define
$$\widetilde{v}=\widetilde{v}(y_1,y_2,\ldots,y_n)=
\cases{v(y_1,y_2,\ldots,y_n), \ \ \ \ \ \ \ \ \ \ \ \ \  \ \ \ y\in
R_+^{n};  \cr v(y_1,y_2,\ldots,y_{n-1},-y_n),\ \ \ \ \ \ y\in
R_{-}^n;}$$
then $\widetilde{v}$ satisfies
\begin{equation}
\left\{\begin{array}{ll}\label{eq18}
 -\Delta \widetilde{v}=\widetilde{v}^p &y\in R^n,\\
\widetilde{v}\geq 0,&y\in R^n,\\
\widetilde{v}(0)=1, \\
\end{array}
\right.
\end{equation}
this also contradicts to Lemma 2.1. $\Box$\\

{\bf Proof of Theorem 1.1}. We argue by contradiction. Suppose that
the conclusion of Theorem 1.1 is false. Then there exists a sequence
$\beta_j\rightarrow 0^+$ as $j\rightarrow\infty$ such that problem
(\ref{eq3}) with $\beta=\beta_j$ has at least one positive solution
$u_{\beta_j}$, by Lemma 2.2 and the standard elliptic estimate,
there exists a positive constant $C$ independent of $j$ such that
$\|u_{\beta_j}\|_{C^{2,\alpha}(\overline\Omega)}\leq C$. Hence, up
to a subsequence, we may assume that $u_{\beta_j}\rightarrow u$ in
$C^2(\Omega)$ as $j\rightarrow\infty$ and $u$ is a solution of the
following problem
\begin{equation}
\left\{\begin{array}{ll}\label{eq19}
 - \Delta u=u^p+f(x),\hspace{0.5cm}
&x\in\Omega,\\
u\geq0,
&x\in\Omega,\\
\frac{\partial u}{\partial \nu}=0,&x\in \partial \Omega.
\end{array}
\right.
\end{equation}
However, it is easy to see that problem (\ref{eq19}) has no solution since
$f(x)\geq 0$ and $f(x)\not\equiv0$. A contradiction. $\Box$

\section{The proof of Theorem 1.2}

\setcounter{section}{3}

\setcounter{equation}{0}

\vskip 0.1in

This section devotes to prove theorem 1.2. Since the proof is
relatively long, we divide it into the following lemmas.\\

{\bf Lemma 3.1.}\ Let $h(x)$ be the solution of problem (\ref{eq7}), and
$\varphi_\beta(x)$ be the solution of the problem
\begin{equation}
\left\{\begin{array}{ll}\label{eq20}
 -\Delta\varphi=1,
&x\in\Omega,\\
\frac{\partial\varphi}{\partial\nu}+\beta\varphi=0,&x\in
\partial\Omega.
\end{array}
\right.
\end{equation}
Then $\varphi_\beta(x)$ is monotonically decreasing with respect to
the parameter $\beta$ and $\lim\limits_{\beta\rightarrow
+\infty}\varphi_\beta(x)=h(x)$.

\vskip 0.1in

{\bf Proof.}\ Let $\beta_1<\beta_2$, $\varphi_{\beta_1}$ and
$\varphi_{\beta_2}$ be the solution of problem (\ref{eq20}) with
$\beta=\beta_1$ and $\beta_2$ respectively. If we set
$w=\varphi_{\beta_1}-\varphi_{\beta_2}$, the $w$ satisfies
\begin{equation}
\left\{\begin{array}{ll}\label{eq21}
 - \Delta w=0,
&x\in\Omega,\\
\frac{\partial w}{\partial \nu}+\beta_1 w\geq 0,&x\in
\partial \Omega.
\end{array}\right.
\end{equation}
It follows from the maximum principle that $w\geq0$ in $\Omega$.
Hence, $\varphi_\beta(x)$ is monotonically decreasing with respect
to $\beta$. Noting that $\varphi_\beta(x)\geq 0$ for any $\beta>0$
and $x\in\Omega$, we conclude that $\varphi_\beta(x)$ converges to
some function. Passing to the limit in Problem (\ref{eq20}) with
$\beta\rightarrow +\infty$, we finally obtain
$$\lim\limits_{\beta\rightarrow +\infty}\varphi_\beta(x)=h(x).$$

\vskip 0.1in

{\bf Lemma 3.2.}\ If $f(x)$ satisfies $(\mathcal{F})$. Then there
exists a positive number $\beta_f^\star$ such that

\vskip 0.1 cm

(i)  Problem (\ref{eq3}) has no solution when
$\beta\in(0,\beta_f^{\star})$.

\vskip 0.1 cm

(ii) If $\beta\geq\beta_f^\star$, then problem (\ref{eq3}) has a minimal
solution $U_\beta(x)$ in the sense that for any solution
$u_\beta(x)$ of problem (\ref{eq3}) we have $U_\beta(x)\leq u_\beta(x)$.
Moreover $U_\beta$ is strictly decreasing with respect to the
parameter $\beta$.

\vskip 0.1in

{\bf Proof.}\ Let $\varphi_\beta$ be the solution of problem (\ref{eq20})
and $\Lambda$ be the constant given in (\ref{eq8}). Setting
$v_\beta=\Lambda\varphi_{\beta}$, it is easy to verify that
$v_\beta(x)$ satisfies
\begin{equation}
\left\{\begin{array}{ll}\label{eq22}
 - \Delta v_{\beta}=-\Lambda\Delta\varphi_{\beta}=\Lambda,
&x\in\Omega,\\
\frac{\partial v_{\beta}}{\partial \nu}+\beta v_{\beta}=0,&x\in\partial \Omega.\\
\end{array}
\right.
\end{equation}
By Lemma 3.1 and the assumption $(\mathcal{F})$, we have
$$\Lambda-\Lambda^{p}\max\limits_{x\in\Omega}\varphi_{\beta}^{p}-f(x)\longrightarrow
\Lambda-\Lambda^pM_h-f(x)>0\ \ \ \ \ \mbox{when}\ \ \beta\longrightarrow +\infty.$$
Consequently, there exists a positive number $\beta^\ast$ such that
 $$\Lambda-\Lambda^p\varphi_{\beta}^{p}-f(x)\geq \Lambda-\Lambda^{p}\max\limits_{x\in\Omega}\varphi_{\beta}^{p}-f(x)>0,$$
for any  $\beta>\beta^\ast$. This and (\ref{eq22}) imply that
$v_\beta=\Lambda\varphi_{\beta}$ is a super-solution of problem
(\ref{eq3}) when $\beta>\beta^\ast$. On the other hand, $0$ is obviously a
sub-solution of problem (\ref{eq3}) and $0<v_\beta(x)$. Hence, the sub and
super-solution method implies that problem (\ref{eq3}) has at least one
solution $u_\beta(x)$ for any $\beta>\beta^\ast$. Moreover,
$0<u_\beta(x)\leq v_\beta(x)$.

Define
$$\beta_f^\star=\inf\{\beta^\ast\in(0,+\infty)\ \mbox{such that problem (1.1) has solution for
any}\ \beta>\beta^\ast\}.$$
It is obvious that problem (\ref{eq3}) has at
least one solution for any $\beta>\beta_f^\star$. Moreover, by a
similar argument to that used in the proof of Lemma 2.2, we can
conclude that for any $\beta\in(\beta_f^\star,\beta_f^\star+1)$,
there exist a positive constant $C$ such that any solution
$u_\beta(x)$ satisfies
$\|u_\beta\|_{C^{2,\alpha}(\overline\Omega)}\leq C$. Hence, up to a
subsequence, we may suppose that $u_\beta(x)$ converges to a
function $u_{\beta_f^\star}(x)$ in $C^2(\Omega)$ as
$\beta\rightarrow\beta_f^\star+0$. Passing to the limit in problem
(\ref{eq3}) with $\beta\rightarrow\beta_f^\star+0$, we know that
$u_{\beta_f^\star}(x)$ is a solution of problem (\ref{eq3}) with
$\beta=\beta_f^\star$. Hence, problem (\ref{eq3}) has at least one
solution for any $\beta\geq\beta_f^\star$. Moreover, we know from
theorem 1.1 that $\beta_f^\star>0$. Next, we prove that problem
(\ref{eq3}) has no solution for any $\beta\in(0,\beta_f^\star)$ by
contradiction. supposing that there exist a number
$\beta_0\in(0,\beta_f^\star)$ such that problem (\ref{eq3}) with
$\beta=\beta_0$ has a solution $u_{\beta_0}$. Then for any
$\beta>\beta_0$, $u_{\beta_0}$ satisfies
\begin{equation}
\left\{\begin{array}{ll}\label{eq23}
 - \Delta u_{\beta_0}(x)=u_{\beta_0}^p(x)+f(x),\hspace{0.5cm}
&x\in\Omega,\\
\frac{\partial u_{\beta_0}(x)}{\partial \nu}+\beta
u_{\beta_0}(x)>\frac{\partial u_{\beta_0}(x)}{\partial \nu}+\beta_0
u_{\beta_0}(x)=0,&x\in
\partial \Omega.
\end{array}
\right.
\end{equation}
This implies that $u_{\beta_0}(x)$ is a super-solution of problem
(\ref{eq3}) for any $\beta>\beta_0$. On the other hand, $0$ is a
sub-solution of problem (\ref{eq3}) and $0<u_{\beta_0}(x)$. Hence, by
the sub-and super-solution method, problem (\ref{eq3}) has at lest one
solution for any $\beta>\beta_0$. This contradicts to the definition
of $\beta_f^\star$. Consequently, problem (\ref{eq3}) has no solution for
any $\beta\in(0,\beta_f^\star)$.

To find a minimal solution of problem (\ref{eq3}), we let $u_0(x)\equiv 0$
and $u_{j+1}(x)$ is the unique solution of the following problem
\begin{equation}
\left\{\begin{array}{ll}\label{eq24}
 - \Delta u_{j+1}=u_j^p+f(x),&x\in\Omega,\\
\frac{\partial u_{j+1}}{\partial \nu}+\beta u_{j+1}=0,&x\in
\partial\Omega.
\end{array}
\right.
\end{equation}
Then it is easy to see that $\{u_j(x)\}$ is a monotonically
increasing sequence. Moreover, if $u_\beta(x)$ is a solution of
problem (\ref{eq3}), we also have $u_j(x)\leq u_\beta(x)$ for any
$x\in\Omega$ and any $j>0$. Hence, $\{u_j(x)\}$ is convergence when
$j$ tends to the infinity. Let

$$U_\beta(x)=\lim\limits_{j\rightarrow +\infty}u_j(x).$$
Then, $U_\beta(x)$ is obviously the minimal solution of problem
(\ref{eq3}).

To prove that $U_\beta(x)$ is monotonically decreasing with respect
to the parameter $\beta$, we assume that $\beta_2>\beta_1$ and
$U_{\beta_2}(x)$ and $U_{\beta_1}(x)$ are the minimal solution of
problem (\ref{eq3}) with $\beta=\beta_2$ and $\beta=\beta_1$ respectively.
Then, a simple computation shows that $U_{\beta_1}(x)$ satisfies
\begin{equation}
\left\{\begin{array}{ll}\label{eq25}
 - \Delta U_{\beta_1}(x)=U_{\beta_1}^p(x)+f(x),&x\in\Omega,\\
\frac{\partial U_{\beta_1}(x)}{\partial \nu}+\beta_2
U_{\beta_1}(x)>\frac{\partial U_{\beta_1}(x)}{\partial \nu}+\beta_1
U_{\beta_1}(x)=0,&x\in
\partial\Omega.
\end{array}
\right.
\end{equation}
Thus, $U_{\beta_1}(x)$ is a super-solution of problem (\ref{eq3}) with
$\beta=\beta_2$. By the definition of $U_{\beta_2}(x)$, we have
$U_{\beta_2}(x)\leq U_{\beta_1}(x)$ for any $x\in\Omega$. Moreover,
by strong maximum principle, we have $U_{\beta_2}(x)<
U_{\beta_1}(x)$ for any $x\in\Omega$. That is, $U_\beta(x)$ is
strictly decreasing with respect to the parameter $\beta$. $\Box$\\

\vskip 0.1in

For any $\beta>\beta_f^\star$, by making use of the minimal solution
of problem (\ref{eq3}), we can decompose any solution $u_\beta(x)$ of
problem (\ref{eq3}) as

$$u_\beta(x)=v(x)+U_\beta(x),$$
where $v(x)$satisfies
\begin{equation}
\left\{\begin{array}{ll}\label{eq26}
 - \Delta v=(U_\beta(x)+v)^p-U_\beta^p(x),&x\in\Omega,\\
\frac{\partial v}{\partial \nu}+\beta v=0,&x\in \partial \Omega.
\end{array}
\right.
\end{equation}

To find another solution of problem (\ref{eq3}), we have only to find a
positive solution of problem (\ref{eq26}). To do this, we represent
$(U_\beta(x)+v)^p-U_\beta^p(x)$ as
$$(U_\beta(x)+v)^p-U_\beta^p(x)=pU_\beta^{p-1}(x)+g(x,v).$$
Then, $g(x,v)$ satisfies

\vskip 0.1in

$(g_1)$\ $\lim\limits_{v\rightarrow 0}\frac{g(x,v)}{v}=0$ uniformly
on $\Omega$;

\vskip 0.1in

$(g_2)$\ $\lim\limits_{v\rightarrow 0}\frac{g(x,v)}{v^p}=0$
uniformly on $\Omega$ with $1<p<\frac{n+2}{n-2}$;

\vskip 0.1in

$(g_3)$\ $g(x,0)=0$ and $g(x,v)>0$ for $v>0)$;

\vskip 0.1in

$(g_4)$\ there exists some $\theta>2$ and $M>0$ such that $0<\theta
G(x,v)\leq vg(x,v)$ for all $v\geq M$ and $x\in\Omega$.

\vskip 0.1in

The following lemma is crucial for finding positive solutions of problem (\ref{eq26}).

\vskip 0.1in

{\bf Lemma 3.3}.\ The first eigenvalue $\lambda_1(\Omega)$ of the
following eigenvalue problem

\begin{equation}
\left\{\begin{array}{ll}\label{eq27}
 -\Delta\phi-pU_\beta^{p-1}\phi=\lambda\phi,&x\in\Omega,\\
\frac{\partial \phi}{\partial \nu}+\beta\phi=0,&x\in\partial\Omega,
\end{array}
\right.
\end{equation}
is positive.

\vskip 0.1in

{\bf Proof.}\ Let $\phi_1(x)$ be the first eigenfunction of problem
(\ref{eq27}). It is well known that $\phi_1(x)$ can be chosen so that
$\phi_1(x)>0$ for any $x\in\Omega$. Let $\beta_1<\beta$,
$U_{\beta_1}(x)$ and $U_\beta(x)$ be the minimal solution of problem
(\ref{eq3}) with parameter $\beta_1$ and $\beta$ in boundary condition
respectively. Then, it follows from lemma 3.2 that

$$U_{\beta_1}(x)>U_{\beta}(x)\ \ \mbox{for}\ \ x\in\Omega.$$
Setting $v(x)=U_{\beta_1}(x)-U_\beta(x)$, then $v(x)>0$ for any
$x\in\Omega$. Furthermore, a simple calculation implies that $v(x)$
satisfies
\begin{equation}
\left\{\begin{array}{ll}\label{eq28}
 -\Delta v=U_{\beta_1}^p-U_\beta^p>pU_\beta^{p-1}v,&x\in\Omega,\\
\frac{\partial v}{\partial \nu}+\beta v=0,&x\in \partial \Omega.
\end{array}
\right.
\end{equation}
From this, we have
\begin{equation}\label{eq29}
\beta\int\limits_{\partial\Omega}v\phi_1ds+\int\limits_{\Omega}\nabla
v\bullet\nabla\phi_1dx-p\int\limits_{\Omega}U_\beta^{p-1}v\phi_1dx>0.
\end{equation}
On the other hand, by the equations satisfied by $\phi_1(x)$, we can
deduce that
\begin{equation}\label{eq30}
\beta\int\limits_{\partial\Omega}v\phi_1ds+\int\limits_{\Omega}\nabla
v\bullet\nabla\phi_1dx-p\int\limits_{\Omega}U_\beta^{p-1}v\phi_1dx=\lambda_1(\Omega)\int\limits_{\Omega}v\phi_1dx.
\end{equation}
Combining (\ref{eq29}) with (\ref{eq30}), we obtain
$$\lambda_1(\Omega)\int\limits_{\Omega}v\phi_1dx>0.$$
Since $\int\limits_{\Omega}v\phi_1dx>0$, we reach
$\lambda_1(\Omega)>0$. $\Box$

\vskip 0.1in

{\bf Proof of Theorem 1.2.}\ The proof of theorem 1.2 (i) is similar
to that of theorem 1.1. The conclusion (ii) of theorem 1.2 comes
from lemma 3.2. To prove (iii) of theorem 1.2, let $H^1(\Omega)$
be the usual Sobolev space endowed with the usual norm
$$\|u\|=\sqrt{\|u\|_{L^2(\Omega)}^2+\|\nabla u\|_{L^2(\Omega)}^2}$$
From lemma 3.3, we can easily see that for $\beta>\beta^*_f$
$$\|u\|_*=\sqrt{\beta\int\limits_{\partial\Omega}v^2ds+\int\limits_{\Omega}|\nabla
v|^2dx-p\int\limits_{\Omega}U_\beta^{p-1}v^2dx}$$ is a norm
equivalent to $\|u\|$. Denote by $H^1_*(\Omega)$ the function space
$H^1(\Omega)$ endowed with the norm $\|\bullet\|_*$. Then
$H^1_*(\Omega)$ is obviously a Banach space. Consider the functional

$$
\begin{array}{ll}
I(v)&=\frac{1}{2}\int\limits_{\Omega}|\nabla
v|^2dx-\frac{1}{2}p\int\limits_{\Omega}U_\beta^{p-1}v^2dx+\frac{\beta}{2}\int\limits_{\partial\Omega}v^2ds
-\int\limits_{\Omega}G(x,v^+)dx\\
&=\frac{1}{2}\|u\|_*^2-\int\limits_{\Omega}G(x,v^+)dx
\end{array}
$$
which is defined on $H^1_*(\Omega)$, where $v^+(x)=\max\{0;v(x)\}$.
By making use of $(g_1) \sim (g_4)$ satisfied by $g(x,v)$, we can
verify that $I(v)$ satisfies all desired conditions in the
Mountain-Pass Theorem with PS-condition (see\cite{ref DengY} for
more details). Thus, Mountain-Pass Theorem with PS-condition
\cite{ref Aarp} assures that problem (\ref{eq26}) has at least one
positive solution for $\beta>\beta^{\star}_f$. Consequently, problem
(\ref{eq3}) has at least two solution for $\beta>\beta^{\star}_f$.
This completes the proof of theorem 1.2. $\Box$

 \section{A priori bound for the global solution of problem (\ref{eq11})}

\setcounter{section}{4}

\setcounter{equation}{0}

\vskip 0.1in

In this section, we derive an a priori bound for global solution of
problem (\ref{eq11}). More precisely, we will prove the following theorem.

\vskip 0.1in

{\bf Theorem 4.1.}\ Assume that $1<p<\frac{n+2}{n-2}$, and $u(x,t)$
is a classical solution of problem (\ref{eq11}) in
$\Omega\times(0,+\infty)$ with $0\leq u_0(x)\in C(\Omega)$. Then
there exists a positive constant $M$ depending only on
$\max\limits_{x\in\Omega}u_0(x)$ and $\sup\limits_{x\in\Omega}f(x)$
such that $u(x,t)\leq M$ for any $(x,t)\in\Omega\times(0,+\infty)$.

When $\beta=+\infty$ and $f(x)\equiv 0$, theorem 4.1 was proved by
Giga in \cite{ref Giga}. Apart from some minor modification, the
proof of theorem 4.1 can be adapted from \cite{ref Giga} line by
line. However, for the reader's convenience and the need of the
sequel paragraph we give a sketch proof here.

\vskip 0.1in

{\bf Lemma 4.1.}\ Let $u(x,t)$ be a nontrivial nonnegative solution
of problem (\ref{eq11}) in $Q=\Omega\times[0,T)$. Suppose that there exists
a positive number $N$ independent of $T$ such that $u(x,t)$
satisfies
\begin{equation}\label{eq31}
\int^T_0\int_\Omega|u_t|^2<N,
\end{equation}
and that for a given $t_0>0$,
$\sup\limits_{x\in\Omega\times[0,T)}u(x,t)$ is attained in
$\Omega\times(t_0,T)$. Then there is a constant $A$ depending only
on $N$ and $t_0$ such that $u(x,t)\leq A$ for any $(x,t)\in Q$.

\vskip 0.1in

{\bf Proof.}\ Suppose the lemma is false. Then there exists a
sequence of solutions $u_k(x,t)$ of problem (\ref{eq11}) with
$T=T_k>0$ and a sequence of points
$(x_k,t_k)\in\Omega\times(t_0,T_k)$ such that
\begin{equation}\label{eq32}
M_k=\sup\limits_{(x,t)\in\Omega\times[0,T_k)}u_k(x,t)=u_k(x_k,t_k)\rightarrow
\infty,\ \mbox{as}\ k\rightarrow\infty.
\end{equation}
and
\begin{equation}\label{eq33}
\int^{T_k}_0\int_\Omega|u_t|^2<N.
\end{equation}
Let $\lambda_k$ be a sequence of positive numbers such that
$$
\lambda_k^{\frac{2}{p-1}}M_k=1.
$$
Then $\lambda_k\rightarrow 0$ as $k\rightarrow\infty$ due to
$M_k\rightarrow+\infty$.

Define a sequence of function $v_k(y,s)$ by
\begin{equation}\label{eq34}
v_k(y,s)=\lambda_k^{\frac{2}{p-1}}u_k(x_k+\lambda_k
y,t_k+\lambda_k^2 s).
\end{equation}

Let $Q_k=\Omega_k\times(-\lambda_k^{-2}t_k,0]$ with
$\Omega_k=\lambda_k^{-1}(\Omega-x_k)=\{y\ \|\
y=\lambda_k^{-1}(x-x_k),\ x\in\Omega\ \}$, and $\partial
Q_k=\partial\Omega_k\times[-\lambda_k^{-2}t_k,0]$. Then, it is easy
to verify that $v_k(y,s)$ satisfies

\begin{equation}
\left\{\begin{array}{ll}\label{eq35}
v_{ks}-\Delta_yv_k=v_k^p+\lambda_k^\frac{2p}{p-1}f(x_k+\lambda_k y)
,&
(y,s)\in Q_k,\\
\frac{\partial v_k}{\partial\nu}+\lambda_k\beta v_k=0,&
(y,s)\in\partial Q_k,\\
v_k(y,s)\leq 1,& (y,s)\in Q_k,\\
v_k(0,0)=1.
\end{array}
\right.
\end{equation}

Let $d_k=dist(x_k,\partial\Omega)$ denote the distance from $x_k$ to
$\partial\Omega$, and $R^n_a=\{y=(y_1,y_2,\cdots,y_n)\in R^n\ \|\
y_n>-a\ \}$. Since $t_k>t_0$, it is easy to see that up to a
subsequence, there hold
$$Q_k\rightarrow R^n\times(-\infty,0]\ \ \ \mbox{or}\ \ \ Q_k\rightarrow
R^n_a\times(-\infty,0]$$ subject to $\lambda_k^{-1}d_k\rightarrow
+\infty$ or $\lambda_k^{-1}d_k\rightarrow a$ respectively. Hence,
our proof are divided into the following two cases.

{\bf Case (i).}\ If $\lambda_k^{-1}d_k\rightarrow +\infty$, then
$Q_k\rightarrow R^n\times(-\infty,0]$. Consequently, for any
parabolic cylinder $Q(r)=\{(y,s)\in R^{n+1}\ \|\ |y|<r,\
s\in(-r^2,0]\ \}$ with radius $r$ in $R^{n+1}$, there exists a
positive integer $K$ such that $Q(r)\subset Q_k$ for all $k>K$.
Since $v_k(y,s)\leq 1$ in $Q_k$, it follows from the parabolic $L^q$
theory \cite{ref  WF} that $v_k(y,s)$ is uniformly bounded in
$W_q^{2,1}(Q(r))$ for any $q>n$. Hence, by the diagonal method, we
can choose a subsequence of $v_k(y,s)$, still denote it by
$v_k(y,s)$, and a function $v\geq 0$ defined on
$R^n\times(-\infty,0]$ such that $v_k(y,s)$ converges uniformly to
$v$ in any parabolic cylinder $Q(r)$. Furthermore, by taking another
subsequence if necessary, we may assume that $v_{ks}(y,s)$ converges
to $v_s$ in $L^2(Q(r))$. Passing to the limit in problem
(\ref{eq35}) with $k\rightarrow +\infty$, and noting that
$\lambda_k^\frac{2p}{p-1}f(x_k+\lambda_k y)\rightarrow 0$ as
$k\rightarrow +\infty$, we see that $v$ solves
\begin{equation}
\left\{\begin{array}{ll}\label{eq36}
v_{s}-\Delta v=v^p,&(y,s)\in R^n\times(-\infty,0],\\
v(y,s)\geq 0,&(y,s)\in R^n\times(-\infty,0],\\
v(0,0)=1.
\end{array}
\right.
\end{equation}
Thanks to (\ref{eq33}), we have $v_s\equiv 0$ in $R^n\times(-\infty,0)$. In
fact, a simple computation shows that
\begin{equation}
\begin{array}{ll}\label{eq37}
\int\limits_{Q(r)}|v_{ks}|^2dyds&=\lambda_k^{\frac{4}{p-1}-n+2}
\int^{t_k}_{t_k-\lambda_k^2r^2}\int_{|x-x_k|<\lambda_kr}|u_t|^2dxdt\\
&\leq\lambda_k^{\frac{4}{p-1}-n+2}\int^{T_k}_0\int_{\Omega}|u_t|^2dxdt\\
&\leq\lambda_k^{\frac{4}{p-1}-n+2}N.
\end{array}
\end{equation}
Since $\frac{4}{p-1}-n+2>0$ due to $p\in(1,\frac{n+2}{n-2})$, we
have
$$\int\limits_{Q(r)}|v_{ks}|^2dyds\rightarrow 0.$$
This yields $v_s=0$ in $Q(r)$ because $v_{ks}$ converges weakly to
$v_s$ in $L^2(Q(r))$ and the norm is lower semi-continuous under
weak convergence. Noting further that $r$ is arbitrary, $v_s$
vanishes identically in $R^n\times(-\infty,0)$. Hence, $v$ is
independent of $s$, and $v\geq 0$ satisfies
\begin{equation}
\left\{\begin{array}{lll}\label{eq38}
 -\Delta v=v^p,&\mbox{in}\ R^n,\\
v(0)=1.\\
\end{array}
\right.
\end{equation}
This contradicts Lemma 2.1.

\vskip 0.1in

 {\bf Case (ii).}\ If $\lambda_k^{-1}d_k\rightarrow a$, then $ Q_k\rightarrow
R^n_a\times(-\infty,0]$. Let $Q_a=R^n_a\times(-\infty,0]$. Then, for
any parabolic cylinder $Q(r)$ in $R^{n+1}$, there exists a positive
integer $K$ such that $Q(r)\cap Q_a\subset Q_k$ for all $k>K$. By
making use of parabolic $L^q$ regularity theory up to the boundary,
as in the case (i) we can choose a subsequence of $v_k$, still
denote it by $v_k$, such that $v_k$ converges uniformly in $Q(r)\cap
Q_a$ to a function $v$ defined on $Q_a$ and such that $v$ satisfies
\begin{equation}
\left\{\begin{array}{ll}\label{eq39}
 -\Delta v=v^p,&y\in \{y_n>-a\},\\
\frac{\partial v}{\partial\nu}=0,&y\in\{y_n=-a\},\\
v(0)=1.
\end{array}
\right.
\end{equation}
Changing coordinates by $y^*=(y_1,y_2,\cdots,y_n+a)$, then
$\hat{v}(y^*)=v(y)$ solves
\begin{equation}
\left\{\begin{array}{ll}\label{eq40}
 -\Delta\hat{v}=\hat{v}^p,&y^*\in R^n_+,\\
\frac{\partial\hat{v}}{\partial\nu}=0,&y^*\in
\partial R^n_+,\\
\hat{v}(0,a)=1.
\end{array}
\right.
\end{equation}

Define
\begin{equation}
\widetilde{v}(y^*)=\left\{\begin{array}{ll}\label{eq41}
\hat{v}(y_1^*,y_2^*,\cdots,y_{n-1}^*,y_n^*),& y^*\in
R_+^{n};\\
\hat{v}(y_1^*,y_2^*,\cdots,y_{n-1}^*,-y_n^*),&y^*\in R_{-}^n;
\end{array}
\right.
\end{equation}
Then $\widetilde{v}$ satisfies

\begin{equation}
\left\{\begin{array}{ll}\label{eq42}
 -\Delta \widetilde{v}=\widetilde{v}^p,& y^*\in R^n,\\
\widetilde{v}(0,a)=1.
\end{array}
\right.
\end{equation}
This also contradicts Lemma 2.1. $\Box$

\vskip 0.1in

{\bf Proof of Theorem 4.1.}\ Assume that
 $u_0\geq 0$ and $u_0(x)\not\equiv 0$. We first note that there are constant $B,\ t^{'}>0$
 depending only on $\sup\limits_{x\in\Omega}u_0(x)$, $|\Omega|$ and $\sup\limits_{x\in\Omega}f(x)$ such that(see \cite{ref WF})

\begin{equation}\label{eq43}
\sup\limits_{0\leq\tau\leq 2t{'}}\sup\limits_{\Omega}u(x,t)\leq B,\
\ \ \int_{\Omega}|\nabla u|^2(x,t^{'})dx\leq B.
\end{equation}
Let
\begin{equation}\label{eq44}
E(t)=E[u]=\frac{1}{2}\int_{\Omega}|\nabla
u|^2dx+\frac{\beta}{2}\int_{\partial\Omega}u^2ds-\frac{1}{p+1}\int_{\Omega}u^{p+1}dx
-\int_{\Omega}fudx
\end{equation}
be the energy associated to problem (\ref{eq11}). Then, it is easy to check
that the following energy identities hold
\begin{equation}\label{eq45}
\frac{d}{dt}\int_{\Omega}|u|^2dx=-4E(t)+\frac{2(p-1)}{p+1}\int_{\Omega}u^{p+1}
-2\int_{\Omega}fudx,
\end{equation}

\begin{equation}\label{eq46}
 \int_{\Omega}|u_t|^2dx=-\frac{d}{dt}E(t).
\end{equation}
By Holder and Young's inequality, we have
\begin{equation}
\begin{array}{ll}\label{eq47}
2\int_{\Omega}fudx &
\leq2\left(\int_{\Omega}u^{p+1}dx\right)^{\frac{1}{p+1}}
\left(\int_{\Omega}f^{\frac{p+1}{p}}dx\right)^{\frac{p}{p+1}}\\
&\leq
\frac{p-1}{p+1}\int_{\Omega}u^{p+1}dx+(\frac{2^{p+1}p}{p-1})^{\frac{1}{p}}\int_{\Omega}f^{\frac{p+1}{p}}dx.
\end{array}
\end{equation}
It follows from (\ref{eq45}) and (\ref{eq47}) that
\begin{equation}\label{eq48}
\frac{d}{dt}\int_{\Omega}|u|^2dx\geq-4E(t)
-(\frac{2^{p+1}p}{p-1})^{\frac{1}{p}}\int_{\Omega}f^{\frac{p+1}{p}}dx+\left(\frac{p-1}{p+1}\right)\int_{\Omega}u^{p+1}dx.
\end{equation}
Identity (\ref{eq46}) says that the energy $E(t)$ should decrease. This
together with (\ref{eq48}) shows that for any $t>0$, we have
\begin{equation}\label{eq49}
E(t)\geq-\frac{1}{4}(\frac{2^{p+1}p}{p-1})^{\frac{1}{p}}\int_{\Omega}f^{\frac{p+1}{p}}dx,
\end{equation}
since otherwise the solution must blow up in finite time.

Integrating (\ref{eq46}) over $(t^{'},T)$ gives
\begin{equation}
\begin{array}{ll}\label{eq50}
\int^T_{t^{'}}\int_{\Omega}|u_t|^2dxdt&=
E(t^{'})-E(T)\\
&\leq
E(t^{'})+\frac{1}{4}(\frac{2^{p+1}p}{p-1})^{\frac{1}{p}}\int_{\Omega}f^{\frac{p+1}{p}}dx\\
&\leq\frac{B}{2}
+\frac{1}{4}(\frac{2^{p+1}p}{p-1})^{\frac{1}{p}}\int_{\Omega}f^{\frac{p+1}{p}}dx.
\end{array}
\end{equation}
Applying Lemma 4.1 with $Q=\Omega\times[t^{'},T)$, we know that any
solution in $\Omega\times[t^{'},T)$ which attains its maximum
outside $\Omega\times[t^{'},2t^{'})$ is bounded from above by a
constant $M$ which is independent of $T$. By (\ref{eq43}) we know that a
solution which take its maximum inside $\Omega\times[t^{'},2t^{'})$
is dominated by $B$. Hence in any case
\begin{equation}\label{eq51}
u(x,t)\leq M+B\ \ \ \mbox{for any}\ \ \
(x,t)\in\Omega\times[t^{'},\infty).
\end{equation}
Combining (\ref{eq43}) with (\ref{eq51}), we reach the desired conclusion of
Theorem 4.1. $\Box$

\section{The proof of Theorem 1.3}

\setcounter{section}{5}

\setcounter{equation}{0}

\vskip 0.1in

This section devotes to prove theorem 1.3. So, we always assume that
$\beta<\beta_f$ with $\beta_f$ being the number determined in
theorem 1.1. If theorem 1.3 is false, then there is a initial data
$u_0(x)$ such that problem (\ref{eq11}) with initial data $u_0(x)$ has a
global solution $u(x,t;u_0)$. For simplicity, we denote $u(x,t;u_0)$
by $u(t)$ if no confusion arise. Let $E(t)$ be the energy associated
to $u(t)$, that is

\begin{equation}\label{eq52}
E(t)=\frac{1}{2}\int_{\Omega}|\nabla
u(t)|^2dx+\frac{\beta}{2}\int_{\partial\Omega}u^2(t)ds-\frac{1}{p+1}\int_{\Omega}u^{p+1}(t)dx
-\int_{\Omega}f(x)u(t)dx.
\end{equation}
A similar argument to that used in the proof of theorem 4.1 shows
that for any $t>0$,
\begin{equation}\label{eq53}
E(t)\geq-\frac{1}{4}(\frac{2^{p+1}p}{p-1})^{\frac{1}{p}}\int_{\Omega}f^{\frac{p+1}{p}}dx
\end{equation}
and
\begin{equation}\label{eq54}
 \int_{\Omega}|u_t(t)|^2dx=-\frac{d}{dt}E(t).
\end{equation}
This yields
\begin{equation}\label{eq55}
\int^{\infty}_{0}\int_{\Omega}|u_t(t)|^2dxdt\leq E(u
_0)+\frac{1}{4}(\frac{2^{p+1}p}{p-1})^{\frac{1}{p}}\int_{\Omega}f^{\frac{p+1}{p}}dx<+\infty.
\end{equation}
Thus, we can pick up a sequence $\{t_j\}$ such that
$t_j\rightarrow\infty$ and $\|u_t(t_j)\|_{L^2(\Omega)}\rightarrow 0$
as $j\rightarrow +\infty$. Moreover, by theorem 4.1, we have
$u(t_j)\leq M$ for some positive constant $M$ independent of $j$.

Multiplying the first equation in problem (\ref{eq11}) by $u(t_j)$ and
integrating over $\Omega$ yields:
\begin{equation}
\begin{array}{ll}\label{eq56}
\beta\int_{\partial\Omega}u^2(t_j)dx+\int_\Omega|\nabla
u(t_j)|^2dx&=\int_\Omega u^{p+1}(t_j)dx+\int_\Omega
f(x)u(t_j)dx-\int_\Omega
u(t_j)u_t(t_j)dx\\
&\leq\int_\Omega u^{p+1}(t_j)dx+\int_\Omega f(x)u(t_j)dx+\int_\Omega
|u(t_j)u_t(t_j)|dx\\
&\leq
M^{p+1}|\Omega|+M\|f\|_{L^1(\Omega)}+M|\Omega|^{\frac{1}{2}}\|u(t_j)\|_{L^2(\Omega)}.
\end{array}
\end{equation}
Taking $\|u_t(t_j)\|_{L^2(\Omega)}\rightarrow 0$ into account, we
arrive
 $$\|u_j\|_{H^1(\Omega)}\leq M^{p+1}|\Omega|+M\|f\|_{L^1(\Omega)}+O(1).$$
From this and the weak compactness of $H^1(\Omega)$, we may assume,
up to a subsequence, that for some function $u(x)\in H^1(\Omega)$
\begin{equation}
\begin{array}{ll}\label{eq57}
u(t_j)\rightharpoonup u(x)\ \ \ \ \mbox{weakly in}\ \ \ \
H^1(\Omega)\\
u(t_j)\rightarrow u(x)\ \ \ \ \mbox{strongly in}\ \ \
L^{p+1}(\Omega).
\end{array}
\end{equation}
For any given $\varphi(x)\in H^1(\Omega)$, multiplying the first
equation in problem (\ref{eq11}) by $\varphi(x)$ and integrating on
$\Omega$ yields
\begin{equation}
\begin{array}{ll}\label{eq58}
\int_\Omega\varphi u_t(t_j)&=\int_\Omega \varphi\Delta
u(t_j)+\int_\Omega u^p(t_j)\varphi+\int_\Omega
f\varphi\\
&=-\beta\int_{\partial\Omega}u(t_j)\varphi-\int_\Omega\nabla
u(t_j)\bullet\nabla\varphi+\int_\Omega u^p(t_j)\varphi+\int_\Omega
f\varphi.
\end{array}
\end{equation}
Thus
\begin{equation}\label{eq59}
\int_\Omega\varphi
u_t(t_j)dx=-\beta\int_{\partial\Omega}u(t_j)\varphi
ds-\int_\Omega\nabla u(t_j)\bullet\nabla\varphi dx+\int_\Omega
u^p(t_j)\varphi dx+\int_\Omega f\varphi dx.
\end{equation}
Taking (\ref{eq57}) and $\|u_t(t_j)\|_{L^2(\Omega)}\rightarrow 0$ into
account, by passing to the limit in (\ref{eq59}) we arrive
\begin{equation}\label{eq60}
\beta\int_{\partial\Omega}u(x)\varphi ds+\int_\Omega\nabla
u(x)\bullet\nabla\varphi dx=\int_\Omega u^p(x)\varphi dx+\int_\Omega
f\varphi dx.
\end{equation}
This implies that $u(x)$ is a weak solution of problem (\ref{eq3}). By
standard regularity theory of elliptic equations, $u(x)$ is also a
classical solution of problem (\ref{eq3}). This contradicts theorem 1.1
and the proof of theorem 1.3 is completed. $\Box$

\section{The Proof of Theorem 1.4}

\setcounter{section}{6}

\setcounter{equation}{0}

\vskip 0.1in

We give a proof of theorem 1.4 in this section. To this end, we need
the following lemmas first.

\vskip 0.1in

{\bf Lemma 6.1.}\ Let $b$ is a positive constant, $p>1$ and
$f(t)=((t+b)^p-b^p)/t$. Then $f(t)$ is monotonically increasing in
$(0,+\infty)$.

\vskip 0.1in

{\bf Proof.}\ An easy computation yields
$$f'(t)=\frac{pt(t+b)^{p-1}-(t+b)^p+b^p}{t^2}.$$
By mean value theorem, there exists $\xi\in(0,t)$ such that
$$(t+b)^p-b^p=pt(\xi+b)^{p-1}.$$
hence
$$f'(t)=\frac{pt[(t+b)^{p-1}-(\xi+b)^{p-1}]}{t^2}.$$
Noting $p>1$, we have $f'(t)\geq 0$ for any $t>0$. $\Box$

\vskip 0.1in

{\bf Lemma 6.2.}\ Assume that $p>1$ and $a>b>0$. If
$f(a)=a^p-pab^{p-1}+pb^p-b^p$, $g(a)=a^p-pa^p+pa^{p-1}b-b^p$, then
$f(a)>0$, $g(a)<0$.

\vskip 0.1in

{\bf Proof.}\ Since $f^{\prime}(a)=p(a^{p-1}-b^{p-1})>0$ and
$f(b)=0$, we have $f(a)>0$. Similarly, we can prove $g(a)<0$. $\Box$

\vskip 0.1in

{\bf Lemma 6.3.}\ Assume that $p>1$ and $a>b>0$. If
$F(\eta)=\eta(a^p-b^p)+b^p-[\eta(a-b)+b]^p$, then $F(\eta)>0$ for
$0<\eta<1$, and $F(\eta)<0$ for $\eta>1$.

\vskip 0.1in

{\bf Proof.}\ An easy computation shows
\begin{equation}
\left\{\begin{array}{ll}\label{eq61}
F^{\prime}(\eta)=a^p-b^p-p(a-b)[\eta(a-b)+b]^{p-1} ,\\
F^{\prime\prime}(\eta)=-p(p-1)(a-b)^2[\eta(a-b)+b]^{p-2}<0.
\end{array}
\right.
\end{equation}
Hence, $F^{\prime}(\eta)$ is strictly decreasing in $[0,+\infty)$.
By lemma 6.2, we have
\begin{equation}
\left\{\begin{array}{ll}\label{eq62}
F^{\prime}(0)=a^p-pab^{p-1}+pb^p-b^p=f(a)>0,\\
F^{\prime}(1)=a^p-pa^p+pa^{p-1}b-b^p=g(a)<0.
\end{array}
\right.
\end{equation}
This implies that there exists a unique point $\xi\in(0,1)$ such
that $F^{\prime}(\xi)=0$. Hence
\begin{equation}
\left\{\begin{array}{lll}\label{eq63}
F^{\prime}(\eta)>0 & \mbox{for} & 0<\eta\leq\xi,\\
F^{\prime}(\eta)<0 & \mbox{for}& \xi<\eta\leq 1.
\end{array}
\right.
\end{equation}
Consequently
\begin{equation}
\left\{\begin{array}{lll}\label{eq64}
F(\eta)>F(0)=0 & \mbox{for}& 0<\eta\leq\xi,\\
F(\eta)>F(1)=0 & \mbox{for}& \xi<\eta<1,\\
F(\eta)<F(1)=0 & \mbox{for} & \eta>1.
\end{array}
\right.
\end{equation}
This is the desired conclusion of lemma 6.3. $\Box$

\vskip 0.1in

The following lemma is important for our proof of theorem 1.4.

\vskip 0.1in

{\bf Lemma 6.4.}\ If $u_1(x)$ and $u_2(x)$ are any two distinct
solutions of problem (\ref{eq3}) which are also different from the minimal
one, then $u_1(x)$ and $u_2(x)$ must intersect somewhere.

\vskip 0.1in

{\bf Proof.}\ Since $u_1(x)$ and $u_2(x)$ are solutions of problem
(\ref{eq3}) which are different from the minimal solution $U_\beta(x)$, we
can decompose $u_1(x)$ and $u_2(x)$ as the following
$$u_1(x)=v_1(x)+U_\beta(x).$$
$$u_2(x)=v_2(x)+U_\beta(x).$$
Moreover, it follows from strong maximum principle that $v_1(x)$ and
$v_2(x)$ are positive in $\overline{\Omega}$. By making use of
equations satisfied by $u_1(x)$ and $u_2(x)$, it is easy to check
that $v_1(x)$ and $v_2(x)$ solve the following problem respectively
\begin{equation}
\left\{\begin{array}{ll}\label{eq65}
 - \Delta v_1=(U_\beta(x)+v_1)^p-U_\beta^p(x) & x\in\Omega\\
\frac{\partial v_1}{\partial \nu}+\beta v_1=0 &x\in \partial \Omega
\end{array}
\right.
\end{equation}
\begin{equation}
\left\{\begin{array}{ll}\label{eq66}
 - \Delta v_2=(U_\beta(x)+v_2)^p-U_\beta^p(x) & x\in\Omega\\
\frac{\partial v_2}{\partial \nu}+\beta v_2=0 &x\in \partial \Omega
\end{array}
\right.
\end{equation}
From these, we can derive
\begin{equation}\label{eq67}
\int\limits_{\Omega}v_2[(U_\beta(x)+v_1)^p-U_\beta^p(x)]dx=\int\limits_{\Omega}v_1[(U_\beta(x)+v_2)^p-U_\beta^p(x)]dx.
\end{equation}
That is
\begin{equation}\label{eq68}
\int\limits_{\Omega}v_1v_2[\frac{(U_\beta(x)+v_1)^p-U_\beta^p(x)}{v_1}-\frac{(U_\beta(x)+v_2)^p-U_\beta^p(x)}{v_2}]dx=0.
\end{equation}
From this and the conclusion of lemma 6.1, we can easily see the
validity of lemma 6.4. $\Box$

\vskip 0.1in

{\bf Proof of Theorem 1.4.}\ The proof of theorem 1.4 (i) is the
same as that of theorem 1.3. To prove theorem 1.4 (ii) and (iii), we
let $u_\beta(x)$ be a solution of problem (\ref{eq3}) which is distinct to
the minimal solution $U_\beta(x)$. Then, by the strong maximum
principle, we have
\begin{equation}\label{eq69}
U_\beta(x)<u_\beta(x)\ \ \ \ \mbox{for any}\ \  x\in\Omega.
\end{equation}
\vskip 0.1in

If the initial datum $0\leq u_0(x)\leq u_\beta(x)$ and
$u_0(x)\not\equiv u_\beta(x)$, then problem (\ref{eq11}) has a
global solution $u(x,t;u_0)$. Moreover, it follows from the strong
comparison principle that $u(x,t;u_0)$ satisfies
\begin{equation}\label{eq70}
u(x,t;u_0)<u_\beta(x)\ \ \ \ \mbox{for any}\ \ (x,t)\in\overline{\Omega}\times(0,+\infty).
\end{equation}
To prove $\lim\limits_{t\rightarrow+\infty}u(x,t;u_0)=U_\beta(x)$,
we may assume, by replacing $u_0(x)$ with $u(x,T;u_0)$ for some
$T>0$ if necessary, that $u_0(x)\leq\eta u_\beta(x)$ for some
$0<\eta<1$. Let
$h_\beta(x)=\eta(u_\beta(x)-U_\beta(x))+U_\beta(x))$. Then, by
making use of lemma 6.3 with $0<\eta<1$, we can verify that
$h_\beta(x)$ satisfies
\begin{equation}
\left\{\begin{array}{ll}\label{eq71}
 - \Delta h_\beta>h_\beta^p+f(x),&x\in\Omega,\\
h_\beta>0,&x\in\Omega,\\
\frac{\partial h_\beta}{\partial \nu}+\beta h_\beta=0,&x\in \partial
\Omega.
\end{array}
\right.
\end{equation}
This implies that $h_\beta(x)$ is a strict super-solution of the
following problem
\begin{equation}
\left\{\begin{array}{ll}\label{eq72}
w_t-\Delta w=w^p+f(x),&(x,t)\in\Omega\times(0,T),\\
\frac{\partial w}{\partial \nu}+\beta w=0,&(x,t)\in \partial \Omega\times[0,T),\\
w(x,0)=h_\beta(x),&x\in\Omega.
\end{array}
\right.
\end{equation}
Let $w(x,t)$ be the solution of problem (\ref{eq72}). Then, $w(x,t)$
exists globally and is strictly decreasing with respect to $t$.
Since $h_\beta(x)>U_\beta(x)$, it follows from the comparison
principle that $w(x,t)>U_\beta(x)$ for any
$(x,t)\in\Omega\times[0,+\infty)$. Hence, we may assume that
$$\lim\limits_{t\rightarrow+\infty}w(x,t)=W(x).$$
Moreover, it is easy to prove that $W(x)$ satisfies
\begin{equation}
\left\{\begin{array}{ll}\label{eq73}
 - \Delta W=W^p+f(x),&x\in\Omega,\\
W>0,&x\in\Omega,\\
\frac{\partial W}{\partial \nu}+\beta W=0,&x\in\partial\Omega.
\end{array}
\right.
\end{equation}
This will imply $W(x)\equiv U_\beta(x)$. Otherwise, on one hand, we
should have $W(x)>U_\beta(x)$, on the other hand, strong comparison
principle should imply $W(x)<u_\beta(x)$ since
$h_\beta(x)<u_\beta(x)$. This contradicts lemma 6.4 since $W(x)$ and
$u_\beta(x)$ are two solutions of problem (\ref{eq3}) which are
distinct to the minimal solution $U_\beta(x)$ and do not intersect
each other. Thus, we have

$$\lim\limits_{t\rightarrow+\infty}w(x,t)=U_\beta(x).$$

Let $v(x,t)$ be the solution of the problem
\begin{equation}
\left\{\begin{array}{ll}\label{eq74}
v_t-\Delta v=v^p+f(x),&(x,t)\in\Omega\times(0,T),\\
\frac{\partial v}{\partial \nu}+\beta v=0,&(x,t)\in \partial\Omega\times[0,T),\\
w(x,0)=0,&x\in\Omega.
\end{array}
\right.
\end{equation}
Then, $v(x,t)$ exists globally and is strictly increasing with
respect to $t$. Moreover, comparison principle implies $v(x,t)\leq
U_\beta(x)$ since $U_\beta(x)$ is a super-solution of problem
(\ref{eq74}). Consequently, $v(x,t)$ converges when $t\rightarrow+\infty$.
It is obvious that
$$\lim\limits_{t\rightarrow+\infty}v(x,t)=U_\beta(x).$$
Since $0\leq u_0(x)\leq\eta u_\beta(x)<h_\beta(x)$, comparison
principle assures that
$$v(x,t)\leq u(x,t;u_0)\leq w(x,t).$$
By applying squeeze principle, we obtain
$$\lim\limits_{t\rightarrow+\infty}u(x,t;u_0)=U_\beta(x).$$
This completes the proof of theorem 1.4 (ii).

Finally, we prove theorem 1.4 (iii) by contradiction. To this end,
we assume that $u_0(x)\geq u_\beta(x)$, $u_0(x)\not\equiv
u_\beta(x)$, and problem (\ref{eq11}) has a global solution $u(x,t,u_0)$.
By the strong comparison principle, we have
\begin{equation}\label{eq75}
u(x,t;u_0)>u_\beta(x)\ \ \ \ \mbox{for any}\ \ \
(x,t)\in\overline{\Omega}\times(0,+\infty).
\end{equation}
Hence, we may assume, by replacing $u_0(x)$ with $u(x,T;u_0)$ for
some $T>0$ if necessary, that $u_0(x)\geq\eta u_\beta(x)$ for some
$\eta>1$. Let $H_\beta(x)=\eta(u_\beta(x)-U_\beta(x))+U_\beta(x))$.
Then, by making use of lemma 6.3 with $\eta>1$, we can verify that
$H_\beta(x)$ satisfies
\begin{equation}
\left\{\begin{array}{ll}\label{eq76}
 - \Delta H_\beta<H_\beta^p+f(x),&x\in\Omega,\\
H_\beta>0,&x\in\Omega,\\
\frac{\partial H_\beta}{\partial\nu}+\beta H_\beta=0,&x\in\partial
\Omega.
\end{array}
\right.
\end{equation}
This implies that $H_\beta(x)$ is a strict sub-solution of the
following problem
\begin{equation}
\left\{\begin{array}{ll}\label{eq77}
G_t-\Delta G=G^p+f(x),&(x,t)\in\Omega\times(0,T),\\
\frac{\partial G}{\partial \nu}+\beta G=0,&(x,t)\in \partial \Omega\times[0,T),\\
G(x,0)=H_\beta(x),&x\in\Omega.
\end{array}
\right.
\end{equation}
Let $G(x,t)$ be the solution of problem (\ref{eq77}). Then, it follows
from the comparison principle that $G(x,t)\leq u(x,t;u_0)$ for any
$(x,t)$ due to $H_\beta(x)<\eta u_\beta(x)\leq u_0(x)$.
Consequently, $G(x,t)$ exists globally and is strictly increasing in
$t$. By theorem 4.1, there exists a positive constant $M$
independent of $t$ such that $G(x,t)\leq M$ for any $(x,t)$. Hence,
$G(x,t)$ converges to some function $g(x)$. Obviously, $g(x)$ is a
solution of problem (\ref{eq3}). Noting that $H_\beta(x)>u_\beta(x)$, it
follows from the strong comparison principle that
$G(x,t)>u_\beta(x)$ for any $(x,t)$. Consequently,
$g(x)>u_\beta(x)$. This contradicts lemma 6.4 since $g(x)$ and
$u_\beta(x)$ are two solutions of problem (\ref{eq3}) which are distinct
to the minimal solution $U_\beta(x)$ and do not intersect each
other.$\Box$

\end{document}